\documentclass[twocolumn,showpacs,preprintnumbers,amsmath,amssymb]{revtex4}
\usepackage{epsfig}
\usepackage{color}
\usepackage{epsf}
\usepackage{float}
\usepackage{hyperref}

\def\reff#1{(\ref{#1})}
\def\emph#1{{\it#1}}

\def\R{{\mathbb R}}
\def\Z{{\mathbb Z}}

\def\G{\bold G}

\def\1{{1\kern-.25em\hbox{\rm I}}}
\def\eu{{1\kern-.25em\hbox{\sm I}}}

\def\R{{\mathbb R}}  
\def\N{{\mathbb N}}  
\def\Z{{\mathbb Z}}  

\def\part#1{
  \shipout\vbox{
    \vbox to 0.35\vsize{}
    \hugetitle #1
    }
  \advancepageno
  \blankpage
}

\def\frac#1#2{{#1\over #2}}

\def\text#1{\quad{\hbox{#1}}\quad}

\def\square{\ifmmode\sqr\else{$\sqr$}\fi}
\def\sqr{\vcenter{
         \hrule height.1mm
         \hbox{\vrule width.1mm height2.2mm\kern2.18mm\vrule width.1mm}
         \hrule height.1mm}}                  

\def\thanks{\noindent{\bf Acknowledgements: }}

\begin{document}
\titlepage

\title{Roughening and inclination of competition interfaces}

\author{Pablo A. Ferrari}
\affiliation{Universidade de S\~ao Paulo} 
\author{James B. Martin}
\affiliation{University of Oxford} 
\author{Leandro P. R. Pimentel}
\affiliation{\'{E}cole Polytechinique F\'{e}d\'{e}rale de Lausanne} 

\date{\today}

\begin{abstract} 
  
We study the \emph{competition interface}
between two clusters growing over a
random vacant sector of the plane in a simple set-up which allows us to perform formal computations and obtain analitical solutions. 
We demonstrate that a phase transition occurs for the asymptotic inclination of this interface when the final macroscopic shape goes from curved to non curved. In the first case it is \emph{random} while in the second one it is deterministic. We also show that the flat case (stationary growth) is a critical point for the fluctuations: for curved and flat final profiles the fluctuations are in the KPZ scale (2/3); for non curve final profile the fluctuations are in the same scale of the fluctuations of the initial conditions, which in our model are Gaussian (1/2).    

\end{abstract}

\pacs{64.60.Ak 64.60.Ht}

\keywords{Suggested keywords}

\maketitle

\paragraph*{\bf{Introduction}}

The behavior of the interface of a growing material has been investigated using
the Eden model \cite{E61}, ballistic deposition and other random systems.  
Typically, the 
growing region converges to an asymptotic deterministic shape and its fluctuations depend on the geometry of the initial condition \cite{kpz86,ps05}. 
A less well studied phenomenon
is the competing growth of two materials. The interface between two growing
clusters (\textit{competition interface}) presents a \emph{random} direction on the
\emph{same} scale as the deterministic shape \cite{dd91,sk95,p04}. In this letter
we describe quite explicitly this phenomenon in a simple model. On gronds of universality, this will provide a guide to understand the interplay between the asymptotics of the competition interface and the final macroscopic shape in models with different growth and competition mechanisms.  

We determine the inclination of the competition interface 
for a growth model called ``last passage percolation'' in a \emph{random} sector of the plane of angle $\theta$. The growth interfaces are mapped into particle configurations of the \textit{totally asymmetric simple exclusion process} in one dimension (TASEP) \cite{r81}.  Under
Euler space-time rescaling, the particle density of the TASEP converges to a
solution of the Burgers equation. This equation has travelling wave solutions
(shocks) corresponding to the case $\theta>180^o$, and rarefaction
fronts corresponding to $\theta<180^o$.  A perturbation at one site of the
initial particle configuration (called a \emph{second class particle}) 
follows a characteristic of the equation or the path of a shock. 
To establish our results, we map the competition
interface linearly onto the path of the second class particle. 

\paragraph*{\bf The growth model}
The random sector is parametrized by the asymptotic slope of its sides. Let
$\lambda\in(0,1]$ and $\rho\in[0,1)$ and define a random path
$\gamma_{0}=(\gamma_{0}(j))_{j\in\Z}\subseteq\Z^{2}$ with $\gamma_{1}(0)=(1,0)$,
$\gamma_0(0)= (1,1)$, $\gamma_{-1}(0)=(0,1)$ as follows.  Starting from $(0,1)$,
walk one unit up with probability $\lambda$ and one unit left with probability
$1-\lambda$, repeatedly, to obtain $\gamma^1_0=(\gamma_{0}(j))_{j<0}$.  Then, starting from
$(1,0)$ walk down with probability $\rho$ and right with probability $1-\rho$ to
get $\gamma^2_0 =(\gamma_{0}(j))_{j>0}$.  $\gamma^{1}_0$ has asymptotic
orientation $(\lambda-1,\lambda)$ while $\gamma^{2}_0$ has asymptotic
orientation $(1-\rho,-\rho)$. Let $C_{0}$ be the sector with boundary
$\gamma_{0}$, containing the first quadrant; its asymptotic angle
$\theta=\theta_{\lambda,\rho}\in[90^{o},270^{o})$ is the angle between
$(\lambda-1,\lambda)$ and $(1-\rho,-\rho)$.  Notice that
$\theta\in[90^{o},180^{o})$ if and only if $\rho<\lambda$.

The path $\gamma_0$ is the growth interface at time 0.  The dynamics are then
defined as follows.  For each $z\in C_0$ and each $t\geq 0$, we have a label
$\sigma_t(z) \in \{0,1,2\}$.  The label is 0 if $z$ is unoccupied at time $t$,
and is 1 or 2 if $z$ belongs to cluster 1 or 2 respectively.  Once occupied, a
site remains occupied and keeps the same value forever.  Initially, set
$\sigma_{0}(z)=1$ for all $z\in\gamma^{1}_0$, $\sigma_{0}(z)=2$ for all
$z\in\gamma^{2}_0$ and $\sigma_{0}(z)=0$ for all $z\in
C_{0}\backslash\gamma_{0}$.  Independently each vacant site $z\in
C_{0}\backslash\gamma_0$ becomes occupied with rate $1$ provided $z-(1,0)$ and
$z-(0,1)$ are occupied.  Let $G(z)$ be the time at which site $z$ becomes
occupied.  At this time $\sigma_t(z)$ assumes the value $\sigma_{t}(\bar{z})$
where $\bar{z}$ is the argument that maximizes $G(z-(1,0))$ and $G(z-(0,1))$.
Thus when a site becomes occupied it joins the cluster of whichever of its two
neighbours (below and to the left) became occupied more recently.  The label of
the site $(1,1)$ may be left ambiguous, but we stipulate that site $(1,2)$
always joins cluster 1, and site $(2,1)$ always joins cluster 2.
 
The process $(\G^{1}_{t},\G^{2}_{t})$, where $\G^{k}_{t}$ is the set of sites
$z\in C_0$ such that $\sigma_{t}(z)=k$, describes the competing spatial growth
model.
\begin{figure}[htb]
\begin{center}
\includegraphics{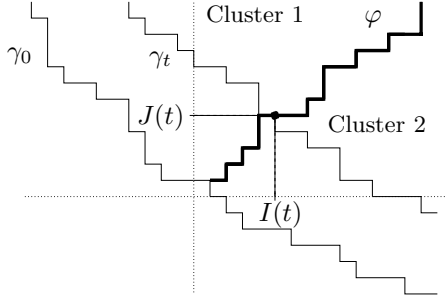}
\end{center}
 \caption{Growth and competition interfaces.}\label{fcompet}
\end{figure}
The growth interface at time $t$ is the polygonal path $\gamma_t$ composed of
sites $z\in C_0$ such that $G(z)\leq t$ and $G(z+(1,1))>t$.  The
\textit{competition interface} $\varphi=(\varphi_{n})_{\N}$ is defined by
$\varphi_{0}=(1,1)$ and, for $n\geq 0$, $\varphi_{n+1}=\varphi_{n} + (1,0)$ if
$\varphi_{n}+(1,1)\in \G_{\infty}^{1}$ and $\varphi_{n+1}=\varphi_{n} + (0,1)$
if $\varphi_{n}+(1,1)\in \G_{\infty}^{2}$. Note that $\varphi$ chooses locally
the shorter step to go up or right, so that $\varphi_{n+1}$ is the argument that
minimizes $\{G(\varphi_n+(1,0)),G(\varphi_n+(0,1))\}$. This competition
interface represents the boundary between those sites which join cluster 1 and
those joining cluster 2 (see Figure \ref{fcompet}). The process
$\psi(t)=(I(t),J(t))$ defined by $\psi(t)=\varphi_n$ for
$t\in[G(\varphi_n),G(\varphi_{n+1}))$ gives the position of the last
intersecting point between the competition interface $\varphi$ and the growth
interface $\gamma_t$.

In \cite{fmp05} we prove that with probability one,
\begin{equation}
  \label{5}
  \lim_{n\to\infty} \frac{\varphi_n}{|\varphi_n|} = e^{i\alpha} 
\end{equation}
where $\alpha\in[0,90^o]$ is given by
\begin{equation}
  \label{p32}
\tan\alpha\;  =\; \left\{\begin{array}{ll} 
 \frac{\lambda\rho}{(1-\lambda)(1-\rho)} &\hbox{ if
 }\rho\ge \lambda\vspace {3mm} 
\\
\bigl(\frac{1-U}{1+U}\bigr)^2&\hbox{ if } \rho< \lambda
\end{array}\right.
\end{equation}
and $U$ is a random variable uniformly distributed in $[1-2\lambda,1-2\rho]$.

\paragraph*{\bf{Simple exclusion and second-class particles}}
The totally asymmetric simple exclusion process $(\eta_{t},\,t\geq 0)$ is a
Markov process in the state space $\{0,1\}^{\Z}$ whose elements are particle
configurations. $\eta_t(j)=1$ indicates a particle at site $j$ at time $t$,
otherwise $\eta_t(j)=0$ (a hole is at site $j$ at time $t$). With rate $1$, if
there is a particle at site $j$, it attempts to jump to site $j+1$; if there is
a hole at $j+1$ the jump occurs, otherwise nothing happens. The basic coupling
between two exclusion processes with initial configurations $\eta_{0}$ and
$\eta'_{0}$ is the joint realization $(\eta_{t},\eta'_{t})$ obtained by using
the same potential jump times at each site for the two different initial
conditions. Let $\eta_{0}$ and $\eta'_{0}$ be configurations of particles
differing only at site $X(0)=0$. With the basic coupling, the configurations at
time $t$ differ only at a single site $X(t)$, 
the position of a so-called \emph{second-class particle}.
Such a particle jumps one step to its right to an empty site with rate $1$, and
jumps backwards one step with rate $1$ when a (first class) particle jumps over
it.

If $\eta_{0}$ is distributed according to the Bernoulli product measure with
density $\lambda$ for $j\leq0$ and $\rho$ for $j> 0$, then the asymptotic
behavior of $X(t)$ shows a phase transition in the line $\lambda=\rho$: with
probability one,
\begin{equation}\label{11} 
  \lim_{t\to\infty} \frac{X(t)}{t} = \left\{\begin{array}{ll} 
1-\rho-\lambda&\hbox{ if } \lambda\le \rho\\ 
U&\hbox{ if }  \lambda> \rho 
\end{array}\right.
\end{equation}
where $U$ is a random variable uniformly distributed in $[1-2\lambda,1-2\rho]$
(\cite{f94,re95,s01} for the deterministic case and \cite{fk95,gm04,fp04,fmp05}
for the random case).

The limits \reff{11} are based on the following hydrodynamic limits.  If
$\eta_{0}$ is distributed with the product measure with densities $\lambda$ and
$\rho$ as before, then the macroscopic density evolution is governed by the
Burgers equation:
\begin{equation}\label{p31}
  \lim_{\epsilon\to 0} \epsilon \sum_{x\in\Z} f(x\epsilon)
  \eta_{t/\epsilon}(x) = \int_\R f(r) u(r,t) dr
\end{equation}
with probability one for all $f:\R\to\R$ with compact support, where $u(r,t)$ is
the solution of the Burgers equation
\[
\frac{\partial u(r,t)}{\partial t}
+\frac{\partial}{\partial r}(u(r,t)(1-u(r,t))) =0,\quad r\in \R,\; t\ge 0
\] 
with initial condition $u(r,0)= \lambda$ for $r\leq 0$ and $\rho$ for $r>0$. If
$\lambda=\rho$ the solution is constant, if $\lambda<\rho$ it is a \emph{shock}:
\begin{equation}\label{d13}
  u(r,t) = \left\{\begin{array}{ll} 
\lambda &\hbox{ if } r\le (1-\lambda-\rho)t\\ 
\rho&\hbox{ if } r> (1-\lambda-\rho)t
\end{array}\right. 
\end{equation}
and it is a \emph{rarefaction front} if $\lambda>\rho$:
\begin{equation}\label{d14}
u(r,t) = \left\{\begin{array}{ll} 
\lambda &\hbox{ if } r\le (1-2\lambda)t\\ 
\frac12- \frac{r}{2(\lambda-\rho)}&\hbox{ if }  (1-2\lambda)t<r\le (1-2\rho)t\\
\rho&\hbox{ if } r> (1-2\rho)t
\end{array}\right. 
\end{equation}
(\cite{r81, re91} for initial product measures and \cite{s98, s01} for initial
measures satisfying \reff{p31} with $t=0$; also \cite{bf94, afs04})

The characteristics $v(a,t)$, corresponding to the Burgers equation and
emanating from $a$, are the solutions of $dv/dt = 1-2u(v,t)$ with $v(0)=a$. The
solutions are constant along the characteristics. When two characteristics
carrying a different solution meet, they give rise to a shock. There is only one
characteristic emanating from locations where the initial data is locally
constant and there are infinitely many characteristics when there is a
decreasing discontinuity. In particular, if the initial condition is
$u(r,0)=\lambda$ for $r<0$ and $u(r,0)=\rho$ for $r\geq 0$, then the
characteristics $v_r(t)$ emanating from the point $v_r(0)=r$ are given by
$v_r(t)=r+(1-2\lambda)t$ if $r< 0$ and $v_r(t)=r+(1-2\rho)t$ if $r>0$.  For
$r=0$ there are two cases. When $\lambda\leq\rho$, the characteristics emanating
from positive sites are slower than those emanating from negative sites. They
collide, giving rise to a shock \reff{d13} traveling at speed $1-\lambda-\rho$.
When $\lambda>\rho$ there are infinitely many characteristics emanating from the
origin: for each $s\in [1-2\lambda,1-2\rho]$ the line $v_0(t) = st$ is a
characteristic emanating from 0. The limits \reff{11} show that the second-class
particle follows the characteristic when there is only one (that is, when
$\lambda=\rho$), that it follows the shock when the initial condition has an
increasing discontinuity and that it chooses uniformly one of the characteristics
emanating from a decreasing discontinuity.

\paragraph*{\bf{Growth and simple exclusion}} 
Rost \cite{r81} relates the simple exclusion process to the growth model as
follows.  Consider initial configurations $\eta_0$ for the exclusion process in
which $\eta_0(0)=0$ and $\eta_0(1)=1$.  Elsewhere let $\eta_0$ be distributed
according to Bernoulli product measure with density $\lambda$ for $j<0$ and
$\rho$ for $j>1$.  Define the initial growth interface $\gamma_0$ by
$\gamma_0(0)=(1,1)$ and $\gamma_{0}(j)-\gamma_{0}(j-1) =
(1-\eta_{0}(j),-\eta_{0}(j))$; then $\gamma_0$ has the same distribution as
before.  Label the particles sequentially from right to left and the holes from
left to right, with the convention that the particle at site $1$ and the hole at
site $0$ are both labeled $1$.  Let $P_j(0)$ and $H_j(0)$, $j\in\Z$ be the
positions of the particles and holes respectively at time 0. The position at
time $t$ of the $j$th particle $P_{j}(t)$ and the $i$th hole $H_{i}(t)$ are
functions of the variables $G(z)$ with $z\in C_{0}\backslash\gamma_{0}$ (defined
earlier for the growth model) by the following rule: at time $G((i,j))$, the $j$th
particle and the $i$th hole interchange positions.  Disregarding labels and
defining $\eta_{t}(P_{j}(t))=1,\,\eta_{t}(H_{j}(t))=0,\, j\in\Z$, the process
$\eta_t$ indeed realizes the exclusion dynamics.  At time $t$ the particle
configuration $\eta_t$ and the growth interface $\gamma_t$ still satisfy the
same relation as $\eta_0$ and $\gamma_0$.  This connection yields the following
\emph{shape theorem} for the growth model.  Almost surely,
\begin{equation}\label{4}
\lim_{t\to\infty}\frac{\gamma_{t}}{t}=\{(r,s)\in\R^{2}\,:\, s= h(r)\}
\end{equation}
where $h(r)=h_{\lambda,\rho}(r)$ is related to the hydrodynamic limit
(\ref{d13},\ref{d14}) by
$h'(r) = {u(r,1)}/{\big(1-u(r,1)\big)}$.

\begin{figure}[htb]
\begin{center}
\includegraphics{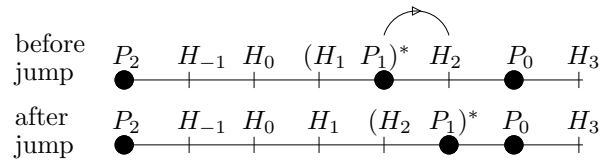}
\end{center}
 \caption{Pair representation of second class particle.}\label{fcompet2}
\end{figure}

\paragraph*{\bf{Second class particles and competition interfaces}}
A key tool in proving (\ref{5},\ref{p32}) is the observation \cite{fp04} that
the process given by the difference of the coordinates of the competition
interface $I(t)-J(t)$ behaves exactly as the second class particle initially put
at the origin. To see this call the particle at site 1 *particle and the hole at
site 0 *hole, and call this couple *pair.  The dynamics of the *pair is the
following: it jumps to the right when the *particle jumps to the right, and it
jumps to the left when a particle jumps from the left onto the *hole.  The *pair
then behaves as a second class particle. The only difference is that it occupies
two sites while the second class particle occupies only one.  The labels of the
*particle and *hole change with time. At time 0 they both have label 1 and the
labels of the *pair are represented by the point $\varphi_0=(1,1)$, the initial
value of the competition interface.  If, say, $G(2,1)<G(1,2)$, then the
*particle jumps over the second hole before the second particle jumps over the
*hole (see Figure \ref{fcompet2}). In this case, the labels of the *pair at time
$G(2,1)$ are $(2,1)$, which is exactly the argument that minimizes
$\{G(2,1),G(1,2)\}$; thus, after the first jump of the *pair, its labels are
given by $\varphi_1$. By recurrence, $\varphi_n$ gives exactly the labels of the
*pair after its $n$th jump. Therefore the labels of the *particle and *hole are
$J(t)$ and $I(t)$, respectively.  In addition, $J(t)-1$ is exactly the number of
jumps that the *pair has made backwards up to time $t$, and $I(t)-1$ is the
number of its jumps forwards.
This shows that if the exclusion and the growth process are realized in the same
space, $X(t) = I(t)-J(t)$.  As a consequence of this and \reff{11} we get the
following behavior for $\psi(t)$ that implies \reff{5} and \reff{p32}. Almost
surely,
\begin{equation}\label{d17}
\lim_{t\to\infty}\frac{\psi(t)}{t}
      = \left\{\begin{array}{ll} 
((1-\rho)(1-\lambda),\lambda\rho)&\hbox{ if } \lambda\le \rho\\ 
\frac14((U+1)^2,(U-1)^2)&\hbox{ if }  \lambda> \rho
\end{array}\right.
\end{equation}
where $U$ is a random variable uniformly distributed in $[1-2\lambda,1-2\rho]$.

To prove \reff{d17} for $\lambda>\rho$ recall that $P_{1}(t)$ is the position of
the $1$st particle at time $t$. Thus $J(t)$ is the number of particles that at
time zero were to the left of $P_{1}(0)=1$ and at time $t$ are to the right of
$X(t)$. Therefore, $J(t)$ is equal to the number of particles between $X(t)$ and
$P_{1}(t)$ at time $t$.  By the law of large numbers, $P_{1}(t)/t$ converges to
$1-\rho$ and, by \reff{11}, $X(t)/t$ converges to $U$. Hence $J(t)/t$ converges
to the integral of the solution of the Burgers equation at time $1$ ($u(r,1)$
given by \reff{d14}) in the interval $(U,1-\rho)$. Taking $f(r) =
\1_{[r\in[U,1-\rho]]}$ in \reff{p31}:
\begin{equation}\label{p30}
  \frac{J(t)}{t}
\;=\;  \frac1t\sum_{j= X(t)}^{P_{1}(t)} \eta_{t}(j) 
\longrightarrow\int_U^{1-\rho} u(r,1) dr = \frac{1}{4}(1-U)^2\nonumber
\end{equation}
Analogously, since $I(t)$ is the number of holes to the right of $H_{1}(0)=0$ at
time zero and to the left of $X(t)$ at time $t$ and $H_{1}(t)/t$ converges to
$-\lambda$ almost surely, we obtain \reff{d17} for $\lambda>\rho$. For
$\lambda\leq\rho$ the same argument works by substituting $U$ above by
$1-\lambda-\rho$, the limit position of the second class particle in this case,
and taking the solution $u(r,1)$ given by \reff{d13}.

\paragraph*{\bf{Fluctuations}}
For $\theta>180^o$ the second class particle has Gaussian fluctuations produced by the initial profile \cite{ff94}. This together with the relation above implies that under a
diffusive scaling $(I(t),J(t))$ converges to a bidimensional Gaussian
distribution with a non-diagonal covariance matrix computed explicitly
\cite{fmp05}.

To understand the fluctuations for $\theta\le 180^o$ we relate the models to a
directed polymer model.  For each $z\in C_{0}\backslash\gamma_0$ let
$w_z=G(z)-\max\{G(z-(1,0)),G(z-(0,1))\}$.  Then $(w_z,\,z\in
C_{0}\backslash\gamma_0)$ is a sequence of i.i.d random variables with an
exponential distribution of mean $1$.  Let $\Pi(z,z')$ be the set of all
directed polymers (or up-right paths) $(z_1,\dots,z_n)$ connecting $z$ to $z'$,
and let $G(z',z)$ be the maximum over all $\pi\in\Pi(z',z)$ of $t(\pi)$, the sum
of $w_z$ along the polymer $\pi$. Each site $z$ has energy $-w_z$, and the
polymer $\pi$ has energy $-t(\pi)$. Thus $-G(z',z)$ is the minimal energy, or
ground state, between $z'$ and $z$. There exists a unique polymer $M(z',z)$ in
$\Pi(z',z)$ that attains the maximum. We say that the semi-infinite polymer
$(z_n)_{\N}$ is \emph{maximizing} if for all $n<m$ we have
$(z_n,\dots,z_m)=M(z_n,z_m)$. Every semi-infinite maximizing polymer $(z_n)_\N$
has an asymptotic inclination $e^{i\alpha}$ \cite{fp04}.  In the competition
model, $G(z)=G(\gamma_0,z)$ and for $k=1,2$, $\G^{k}_\infty$ is the set of sites
$z$ such that $M(\gamma_0,z)$ originates from $\gamma^{k}_0$.  Denote by $\xi$ the
roughening exponent of semi-infinite maximizing polymers.

For $\theta=180^o$ ($\lambda=\rho$) the process is stationary and the connection
is explicit. Running the process forward and backward we extend $G(z)$ to all
$z\in\Z^2$; $G^+=(G(z),\,z\in\Z^2)$ and $G^-=(-G(-z),\,z\in\Z^2)$ are
indentically distributed. We define the \emph{forward} competition interface
starting at $z$, $\varphi^{z}=(\varphi^{z}_n)_\N$, by setting $\varphi^{z}_0=z$
and putting $\varphi^{z}_{n+1}$ equal to the argument of the minimum between
$G(\varphi^{z}_n +(1,0))$ and $G(\varphi^{z}_n +(0,1))$, and the \emph{backward}
semi-infinite polymer starting at $z$, $M^{z}=(M^{z}_n)_\N$, by setting
$M^{z}_0=z$ and putting $M^{z}_{n+1}$ equal to the argument of the maximum
between $G(M^{z}_n -(1,0))$ and $G(M^{z}_n -(0,1))$. Note that
$\varphi=\varphi^{(1,1)}$ and that $M^{z}$ is a semi-infinite maximizing
polymer.  Together with the duality relation $\varphi^{z}(G^+)=M^{z}(G^-)$, this
shows that the forward competition interface has the same law as the backward
semi-infinite maximizing polymer and, in particular, they have the same
fluctuations, so that $\chi=\xi$.

For $\theta<180^o$ ($\lambda>\rho$) the competition interface $\varphi$ is
enclosed by two semi-infinite maximizing polymers $M^1$ and $M^2$ starting from
$\gamma^1_0$ and $\gamma^2_0$, respectively, and with the same inclination
\cite{fmp05}.  Therefore $\chi\leq\xi$ in this case.

\paragraph*{\bf Conclusions} 

The connections studied above between the competition interface, the
second class particle and maximal polymers fit into the interplay between
the fluctuation statistics and the global geometry of the growth interface
developed by Prahofer and Spohn \cite{ps05}. If the final macroscopic
profile is curved then the competition interface follows a random direction
(characteristic) intersecting the final surface at a point with non-zero
curvature. In this case we have the KPZ scaling and the competition
interface gets the transversal fluctuations, indicating the exponent
$\chi=2/3$. If the macroscopic profile is not curved we have two different
situations. In the flat case (stationary growth) the competition interface
also follows the characteristics of the associated hydrodynamic PDE and we
still have the KPZ scaling. In the shock case the competition interface
gets the longitudinal fluctuations which, in this case, are produced by
the Gaussian fluctuations ($\chi=1/2$) of the initial profile. On microscopic grounds one might suggest different rules for growth and competition. By universality we expect that from the knowledge of the curvature of the final macroscopic shape one can infer the asymptotics of the competition interface. This fits with the exponents founded by Derrida and Dickman \cite{dd91} in the Eden context since, in this case, the macroscopic profile is curved for angles $\theta>180^o$ (we notice that in their simulations they have considered periodic initial conditions and so the longitudinal fluctuations in the shock direction are governed by the exponent $1/3$ \cite{ps05}).       
         
\paragraph*{\bf{Acknowledgments}} We thank R. Dickman for calling our attention
to this problem and a referee for a careful reading and useful comments about a previous version of this paper.

\end{document}